\documentclass[a4paper,landscape]{article}

\usepackage{a4} 
\usepackage{amsfonts}
\usepackage{amssymb}
\usepackage{amsmath}
\usepackage{amscd}
\usepackage{amsthm}
\usepackage{epsfig}
\usepackage[all]{xy}
\usepackage{color}
\usepackage{bbm}
\usepackage{graphicx}
\usepackage{psfrag}
\usepackage{hyperref}
\usepackage{enumerate}

\newcommand{\R}{\mathbb{R}}

\newcommand{\bprf}{\begin{proof}[Proof]}
\newcommand{\eprf}{\end{proof}}

\newcommand{\ra}{\rightarrow}

\newcommand{\conv}{{\rm conv}}
\newcommand{\cone}{{\rm cone}}

\graphicspath{{pictures/}}   

\newtheoremstyle{plainNoItalics}{}{}{\normalfont}{}{\bfseries}{.}{ }{}

\newtheorem{theorem}{Theorem}

\newtheorem{lemma}[theorem]{Lemma}

\newtheorem*{theorem*}{Theorem}
\newtheorem*{proposition*}{Proposition}
\newtheorem*{lemma*}{Lemma}
\newtheorem*{corollary*}{Corollary}

\theoremstyle{plainNoItalics}

\newtheorem*{definition*}{Definition}
\newtheorem{remark}[theorem]{Remark}

\newtheorem*{remark*}{Remark}

\newtheorem*{observation*}{Observation}
\newtheorem*{example*}{Example}
\newtheorem*{assumption*}{Assumption}
\newtheorem*{problem}{Problem}

\begin{document}
\title{On the number of Birch partitions}
\author{{\Large Stephan Hell}}
\date{Institut f\"ur Mathematik, MA 6--2,
TU Berlin,\\ D--10623 Berlin, Germany,
hell@math.tu-berlin.de}

\maketitle

\begin{abstract}
  Birch and Tverberg partitions are closely related concepts from
 discrete geometry. We show
  two properties for the number of Birch partitions: Evenness, and a
  lower bound. This implies the first non-trivial lower bound for the
  number of Tverberg partitions that holds for arbitrary $q$, where $q$ is the number of partition blocks. The
  proofs are based on direct arguments, and do not use the equivariant
  method from topological combinatorics. 
\end{abstract}
\section{Introduction}\label{sec-intro}
Our starting point is the following theorem due to B.~J.~Birch~\cite{birch59} from 1959.
\begin{theorem}\label{theorem-birch}Given $3N$ points in $\R^2$, we
  can divide them into $N$ triads such that their convex hulls contain
  a common point.
\end{theorem}
The proof of Theorem~\ref{theorem-birch} is based on a lemma on
partitioning a general measure which is due to Richard Rado, nowadays
known as the {\it center point theorem}. See e.g.~Matou\v{s}ek's textbook~\cite{matousek02:_lectur_discr_geomet}, 
or Tverberg and Vre\'cica~\cite{tverberg93} for more details.

Theorem~\ref{theorem-birch} led us to the following definition, 
see also Tverberg and Vre\'cica~\cite{tverberg93}.
\begin{definition*}\label{def-birch}
  Let $X$ be a set of $k(d+1)$ points in $\R^d$ for some $k\geq 1$. A
  point $p\in\R^d$ is a {\it Birch point} of $X$ if there is a
  partition of $X$ into $k$ subsets of size $d+1$, each containing $p$
  in its convex hull. The partition of $X$ is a {\it Birch partition
    for $p$}. For fixed $p\in\R^d$, let $B_p(X)$ be the number of
  unordered Birch partitions for $p$.
\end{definition*}
From now on, we fix $p$ to be the origin, and we write {\it Birch
  partition} instead of {\it Birch partition for the origin} for
short. A set of points in $\R^d$ is {\it in general position} if
no $k+2$ points are on a common $k$-dimensional affine
subspace. A set $X$ of points in $\R^d$ is {\it in general position with respect to a point $p$} if $X\cup\{p\}$ is in general position.

Our first main result is the following theorem on the number of Birch
partitions.
\begin{theorem}\label{thm-number-birch-part}Let $d\geq 1$ and $k \geq 2$ be
  integers, and $X$ be a set of $k(d+1)$ points in $\R^d$ in general
  position with respect to the origin $0$. Then the following
  properties hold for $B_0(X)$:
\begin{enumerate}[\rm i)]
\item\label{it-birch-even}$B_0(X)$ is even.
\item\label{it-birch-lower}$B_0(X)>0\,\,\Longrightarrow\,\,B_0(X)\geq k!$
\end{enumerate}
\end{theorem}
If the origin is not in the convex hull of $X$, then one has $B_0(X)=0$
which is even. If there is a Birch partition then the lower bound given
in Property~\ref{it-birch-lower}) is tight. Based on computer experiments,
we moreover conjecture:
\begin{eqnarray}\label{bp-conj}B_0(X)\leq (k!)^d.
\end{eqnarray}

B.~J.~Birch proved Theorem~\ref{theorem-birch} to obtain the following
statement for $d=2$. Helge Tverberg then settled the problem for arbitrary
dimension $d$ in 1966.
\begin{theorem}[Tverberg's theorem]\label{thm-tverberg}
  Let $d$ and $q$ be integers. Any $(q-1)(d+1)+1$ points in $\R^d$ can
  be partitioned into $q$ subsets such that their convex hulls have a
  point in common.
\end{theorem}
Partitions as in Theorem~\ref{thm-tverberg} are {\it Tverberg
partitions (into $q$ blocks)}. 
From now on, we implicitly assume that Tverberg partitions
are partitions into $q$ blocks given a set of $(q-1)(d+1)+1$ 
points in $\R^d$. The point in common is a {\it Tverberg point}.

A wave of excitement started in 1981 when B\'ar\'any et al.~\cite{bss81:_tverb} were able to
prove a more general topological version known as the {\it Topological
Tverberg Theorem} when $q$ is a prime number using Borsuk-Ulam's
theorem from algebraic topology. This has then been extended to
prime powers $q$ by many authors, e.~g.~\"Ozaydin~\cite{oezaydin87:_equiv}, 
Volovikov~\cite{volovikov96:_tverb_theor}, 
Sarkaria~\cite{sarkaria00:_tverb_borsuk}. The general
case for arbitrary $q$ is still open; see Matou\v{s}ek's 
textbook~\cite{matou03:_using_borsuk_ulam} for more background.
 
The number of Tverberg partitions has
been studied by Vu\'ci\'c and \v{Z}ivaljevi\'c~\cite{vz93:_notes_sierk}, 
and Hell~\cite{hell07:_tverb}. Using the
equivariant method from topological combinatorics they have obtained:
\begin{theorem}\label{theorem-top-lower-tp}Let $q=p^r$ be a prime 
  power and $d\geq 1$. For any continuous map \mbox{$f:\sigma^N\ra
    \R^d$}, where $N=(d+1)(q-1)$, the number of unordered $q$-tuples
  $\{F_1,F_2,\ldots,F_q\}$ of disjoint faces of the $N$-simplex
  with\linebreak $\bigcap_{i=1}^{q}f(\| F_i \|)\not=\emptyset$ is at
  least
\[\frac{1}{(q-1)!}\cdot\left(\frac{q}{r+1}\right)^{\lceil\frac{N}{2}\rceil}.\]
\end{theorem}
Restricting $f$ to an affine map, unordered $q$-tuples as in 
Theorem~\ref{theorem-top-lower-tp}
are in bijection with Tverberg partitions of the set 
$\{f(v_i)\,|\,v_0,v_1,\ldots,v_N \text{ vertices of } \sigma^N\}$
of $N+1=(d+1)(q-1)+1$ many points in $\R^d$.

Using Theorem~\ref{thm-number-birch-part}, we obtain our second main 
result: The first
non-trivial lower bound for the number of Tverberg partitions
that holds for arbitrary $q$.
\begin{theorem}\label{thm-lower-aff-arb}Let $X$ be a set of
  $(d+1)(q-1)+1$ points in general position in $\R^d$. Then
  the following properties hold for the number $T(X)$ of Tverberg
  partitions: \begin{enumerate}[\rm i)] 
  \item $T(X)$ is even for $q>d+1$.  
  \item\label{item-tp-lower} $T(X)\geq (q-d)!$ 
  \end{enumerate}
\end{theorem}
Property~\ref{item-tp-lower}) improves the result of Theorem
\ref{theorem-top-lower-tp} for $d=2$ and $q\geq 7$. Sierksma
conjectured in 1979 that $T(X)$ is bounded from below by
$((q-1)!)^d$. Combining Theorem~\ref{thm-lower-aff-arb} 
and methods from topological combinatorics, we have been
able to confirm this conjecture for $d=2$ and $q=3$ 
in Hell~\cite{hell08:_tverberg_constraints}, see also 
Hell~\cite{hell06:_tverb_fract_helly}.\\

In Section~\ref{sec-number-bp}, we prove Theorem
\ref{thm-number-birch-part}. Section~\ref{sec-new-lower-tp} comes with
a proof of Theorem~\ref{thm-lower-aff-arb}.  

\section{On the number of Birch partitions}
\label{sec-number-bp}
\begin{figure}[h]
  \centering
  \includegraphics{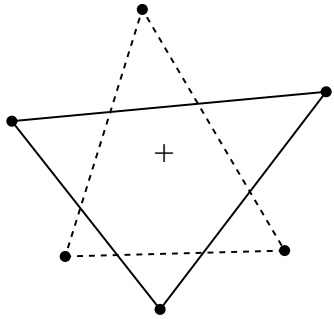}
  \caption{A Birch partition for $6$ points in the plane. }
  \label{fig-birch}
\end{figure}

Figure~\ref{fig-birch} shows a Birch partition for the origin 
denoted as $+$. Each triangle
corresponds to a partition block. There is another way to obtain a
Birch partition for the origin in this example. 

For $d=1$, a Birch partition of a set $X$ of $2k$ points corresponds to 
$k$ intervals containing $0$. Therefore $k$ points of X are in $\R^+$, and
$k$ many in $\R^-$. It is easy to check that there are exactly $k!$
ways to obtain a Birch partition. Hence we have settled Theorem~\ref{thm-number-birch-part}
for $d=1$.\\

We now prove Theorem~\ref{thm-number-birch-part} for $d\geq 2$ in two
steps: We first prove Property~\ref{it-birch-even}), then we prove that
Property~\ref{it-birch-even}) implies Property~\ref{it-birch-lower}).

In our proof, we make use of the 
following basic lemma; see e.~g.~B\'ar\'any and Matou\v{s}ek~\cite{barany07:_quadr}, 
or Deza et al.~\cite{deza06:_colour} for a proof.
\begin{lemma}\label{lemma-convex-cone}
  If $X\subset\R^d$ is a set of points in
  general position with respect to the origin~$0$ and $p\in X$, 
  then $0\in\conv (X)$ if and only  if $-p\in\cone(X\setminus\{p\})$.
\end{lemma}

\begin{figure}[h]
  \centering
  \includegraphics{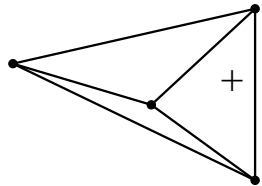}
  \caption{Four points that form two triangles containing the origin. }
  \label{fig-3simpl}
\end{figure}

The following lemma is an easy consequence of Lemma~\ref{lemma-convex-cone}.
\begin{lemma}\label{lem-even-convex}Let $X$ be a set of $d+2$
  points in $\R^d$ that is in general position with respect to the
  origin. Then the number of $d$-simplices with vertices in $X$ that
  contain the origin is even. In fact, this number is either $0$, or
  $2$.
\end{lemma}
See Figure~\ref{fig-3simpl} for a configuration of four points 
in dimension $d=2$ such that two triangles contain the origin $+$. 

\bprf (of Theorem~\ref{thm-number-birch-part}) We first prove
Property~\ref{it-birch-even}) for arbitrary $d\geq 2$, by induction on $k\geq
2$. The base case $k=2$ is the key
part.\\

$k=2$: 
If all rays of the $2d+2$ points of $X$ intersect $S^{d-1}\subset\R^d$ 
close to the north pole then $B_0(X)=0$, as $0\not\in\conv(X)$. We move one
point $p$ of $X$ at a time while all other points remain fixed. The point $p$
can be moved on its ray without changing $B_0(X)$.
Instead of following $p$, we look at its antipode $-p$ as for any
$d$-element subset $S$ of $X\setminus\{p\}$ one has due to Lemma
\ref{lemma-convex-cone}:
\[ 0\in\conv(S\cup\{p\}) \text{\;\; iff \;\;} -p\in\cone(S).\] Every
$d$-element subset of $X\setminus\{p\}$ defines a cone, and these
cones define a decomposition of the sphere $S^{d-1}\subset\R^d$ into cells. The
boundary of a cell is defined by hyperplanes spanned by
$(d-1)$-element subsets of $X\setminus\{p\}$ and the origin. At some
point we are forced to move $-p$ transversally from one side of a
boundary hyperplane defined by a $(d-1)$-element subset $T$ to the
other side. When $-p$ crosses such a hyperplane then $B_0(X)$ might
change. We show in the case distinction below that for every change
the parity of $B_0(X)$ does not change. The number $B_0(X)$ is thus
even as we can move every point of $X$ to its position while fixing
all other points. The cell decomposition during this process is nice:
We can move $-p$ to every position on the sphere while crossing
hyperplanes in a transversal way.

Let us first look at the set of all $d$-simplices $S$ spanned by $d+1$
points from $X$ that contain the origin. If $-p$ crosses the
hyperplane through $T$ transversally, this set might change. For this,
put $\tilde{T}=T\cup\{p\}$.  For all simplices that do not contain
$\tilde{T}$ as a face nothing changes.  If $S$ is of the form
$\tilde{T}\cup\{x\}$ for some $x\in X\setminus\tilde{T}$, then this
property switches:
\begin{quote}$0\in\conv (S)$
before the crossing iff $0\not\in\conv (S)$ afterwards.
\end{quote}

A Birch partition consists of a $d$-simplex $S$ and its complement
$\bar{S}$ in $X$ -- which is again a $d$-simplex -- such that both
contain the origin.  The change of $B_0(X)$ coming from the crossing
of $-p$ can thus only be affected by partitions that contain
$\tilde{T}$ as a face of $S$, or of $\bar{S}$.

Case 1: The complements of all simplices using $\tilde{T}$ do not
contain the origin. $B_0(X)$ does not change as the set of all Birch
partition remains the same.

Case 2: Assume that $\tilde{T}$ is not part of a $d$-simplex $S$ such
that $\{S,\bar{S}\}$ is a Birch partition, and that after the crossing
of $-p$ a Birch partition comes up. We show that Birch partitions come
up in pairs. 

Suppose there is a new Birch partition of the form
$S=\tilde{T}\cup\{x_1\}$ together with its complement $\bar{S}$. Due
to Lemma~\ref{lem-even-convex} there is exactly two $d$-simplices in
$\bar{S}\cup\{x_1\}$ such that both contain the origin. One of them is
$\bar{S}$, let $S^*$ be the other. By assumption $0\not\in\bar{S^*}$
before the crossing of $-p$. In fact, $\bar{S^*}=\tilde{T}\cup\{x_2\}$
for some $x_2$. The set $\{\bar{S^*},S^*\}$ is thus our second Birch
partition as $0\in\conv(\bar{S^*})$ afterwards.\\
Suppose there are three Birch partitions of the form
$S_1=\tilde{T}\cup\{x_1\}$, $S_2=\tilde{T}\cup\{x_2\}$, and
$S_3=\tilde{T}\cup\{x_3\}$, with $x_1,x_2,x_3\in X\setminus\tilde{T}$,
together with their complements.  This can not happen: One has
$0\in\bar{S_i}$ for $i=1,2,3$, and $|\bigcup_{i=1}^3\bar{S_i}|=d+2$.
This contradicts Lemma~\ref{lem-even-convex}. Hence the two new Birch
partitions are of the form $\tilde{T}\cup\{x_1\}$
resp.~$\tilde{T}\cup\{x_2\}$, with $x_1,x_2\in X\setminus\tilde{T}$,
plus their complements.

Case 3: This is the inverse case of Case 2. Assume that there are
exactly two Birch partitions of the form $\tilde{T}\cup\{x_1\}$
resp.~$\tilde{T}\cup\{x_2\}$, with $x_1,x_2\in X\setminus\tilde{T}$,
plus their complements before the crossing. Both of them vanish
after crossing of $-p$. New Birch partitions do not come up as
for this we needed another $\tilde{T}\cup\{x_3\}$ such that
its complement contains the origin. This cannot exist due
to Lemma~\ref{lem-even-convex}.

Case 4: Assume there is exactly one Birch partition of the form
$S=\tilde{T}\cup\{x\}$, with $x\in X\setminus\tilde{T}$, together with
its complement before the crossing. This Birch partition vanishes, and
a new one comes up. 

One has $0\not\in S$ after the crossing of $-p$ so that
$\{S,\bar{S}\}$ vanishes. As in Case 2, there are exactly two
$d$-simplices in $\bar{S}\cup\{x\}$ such that each contains the
origin.  One of them is $\bar{S}$, let $S^*$ be the other. By
assumption $0\not\in\bar{S^*}$ before the crossing of $-p$. In fact,
$\bar{S^*}=\tilde{T}\cup\{x'\}$ for some $x'$. The set
$\{\bar{S^*},S^*\}$ is thus the new Birch partition as
$0\in\conv(\bar{S^*})$ afterwards.\\

Let now $k\geq 3$, and let $p$ be a point in $X$. Let $F_1^{(1)},F_1^{(2)},\dots, F_1^{(l)}$
be all $d$-simplices containing $p$ that can be completed to a Birch
partition of the origin into $k$ subsets. For every $F_i$, omitting
$F_i$ leads to a Birch partition into $k-1$ subsets. By induction
hypothesis, there is an even number of Birch partitions into $k-1$
subsets for the restriction of every $F_i$. \\

Now we assume Property~\ref{it-birch-even}), and derive Property~\ref{it-birch-lower}) 
by induction on $k\geq 2$. The case $k=2$ is due
to Property~\ref{it-birch-even}): $B_0(X)$ is even, so
\[ B_0(X)>0\,\,\Longrightarrow\,\,B_0(X)\geq 2=k!\] Let $k\geq 3$ and
$B_0(X)>0$. Then there is a Birch partition $F_1,F_2,\ldots ,F_k$. If
we take any $k-1$ of the $F_i$, they form again a Birch partition. By
induction hypothesis, the union of $k-1$ many $F_i$ has at least
$(k-1)!$ Birch partitions. In particular, there are $(k-1)!$ many Birch
partitions of $X$ into $k$ subsets that start with $F_1$.  Let $p$ be an
element of $F_1$.

\begin{figure}[h]
  \centering
  \includegraphics{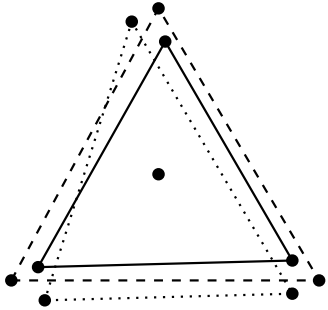}
  \caption{A planar configuration with
    $36=(3!)^2$ Birch resp.~Tverberg partitions.}
  \label{fig-sierksma}
\end{figure}

For every pair $F_1,F_i$, for $i\in\{2,3,\ldots k\}$, one has again
$B_0(F_1\cup F_i)>0$ and so there is a second Birch partition
$\tilde{F}_1^i,\tilde{F}_i$ of $F_1\cup F_i$. Assume without loss of
generality $p\in\tilde{F}_1^i$. The $k$ sets
$F_1,\tilde{F}_1^2,\tilde{F}_1^3,\ldots,\tilde{F}_1^k$ are pairwise
distinct by construction. Every one of them contributes $(k-1)!$ many
Birch partitions of $X$ by induction hypothesis.\eprf
\begin{remark}\label{rem--bp-lower-bound}
  In the induction of our second step, we
  didn't make use of convexity. The key is the base case $k=2$:
\[B_0(X)>0\,\,\Longrightarrow\,\,B_0(X)\geq 2.\]
\end{remark}
\begin{remark*}\label{rem-sierksma-upper-birch}
  Sierkma's configuration shown for $d=2$ and $q=4$ 
  in Figure~\ref{fig-sierksma} attains the
  conjectured upper bound~(\ref{bp-conj}) for $B_p(X)$. Hence it would be
  \emph{maximal} for the number of Birch partitions. At the same time,
  Sierksma conjectured it to be \emph{minimal} for the number of Tverberg
  partitions.
\end{remark*}

\section[On the number of Tverberg partitions]
{On the number of Tverberg partitions}
\label{sec-new-lower-tp}
In this section, we prove Theorem~\ref{thm-lower-aff-arb}. The proof 
is based on the fact that Birch partitions come up while 
studying Tverberg partitions.

Figure~\ref{fig-sierksma} shows a set $X$ of 
$(d+1)(q-1)+1=10$ points in the plane
for $q=4$. A Tverberg partition can be read off as follows:
Each triangle corresponds to a partition block. The point
in the center is the forth block, and at the same time a Tverberg point. 
 
In our proof, we need the following
reformulation of Lemma~2.7 from Sch\"oneborn and 
Ziegler~\cite{schoeneborn05:_topol_tverb_theor}. 
\begin{lemma}\label{obs-gen-pos-tp}
  Let $X$ be set of $(d+1)(q-1)+1$ in
  general position in $\R^d$. Then a Tverberg partition consists of:\\
  \textbullet\,\, Type~I: One vertex $v$, and $(q-1)$ many
  $d$-simplices containing $v$. \\
  \textbullet\,\, Type~II: $k$ intersecting simplices of dimension
  less than $d$, and $(q-k)$ $d$-simplices containing the
  intersection point for some $1<k\leq \min\{d,q\}$.
\end{lemma}
For $d=2$, a type~II partition consist of two intersecting segments,
and $q-2$ many triangles containing their intersection point.

\bprf (of Theorem~\ref{thm-lower-aff-arb}) Tverberg's Theorem
\ref{thm-tverberg} implies the existence of a Tverberg partition
together with a Tverberg point $p$.  The set $X$ is in general
position such that the partition is either of type~I, or type~II.

For type~I, $q-1$
disjoint $d$-simplices contain a point $p$ of $X$.  The $q-1$ disjoint
$d$-simplices make up a Birch partition for $p$. Theorem
\ref{thm-number-birch-part} implies that there are at least $(q-1)!$
many Birch partitions of $p$.  Hence there are at least $(q-1)!$ many
Tverberg partitions.

For type~II, the Tverberg point $p$ is the
intersection of the convex hull of $k\leq d$ many sets of cardinality
at most $d$. The remaining points are partitioned into $q-k$ many
$d$-simplices containing $p$. For $q>d+1$, this makes up a Birch
partition for $p$ into $q-k\geq 2$ sets.  Again by Theorem
\ref{thm-number-birch-part} there are at least $(q-k)!$ Tverberg
partitions.

Properties i) and ii) follow from the corresponding
results on the number of Birch partitions from Theorem
\ref{thm-number-birch-part}. For $q>d+1$, both types of Tverberg
partitions correspond bijectively to Birch partitions so that the
number of Tverberg partitions is even. As we can not predict the type
of the Tverberg partition, the lower bound is equal to $(q-d)!$.
\eprf
\begin{remark*}\label{rem-prf-lower}
\begin{enumerate}
\item Our proof shows a bit more than a lower bound of $(q-d)!$. If we knew
what type of Tverberg partition showed up, then we would obtain $(q-k)!$
for some $k\in\{1,2,\ldots,d\}$. If there is a Tverberg partition of
type~I then the lower bound equals $(q-1)!$.
\item In Hell~\cite{hell08:_tverberg_constraints}, we improve the result
of Theorem~\ref{thm-lower-aff-arb} by proving a lower bound
for the number of Tverberg points, and by using Tverberg's theorem
with constraints.
\end{enumerate}
\end{remark*}

\section{Further directions}
\begin{figure}[h]
    \centering 
    \includegraphics[width=4cm]{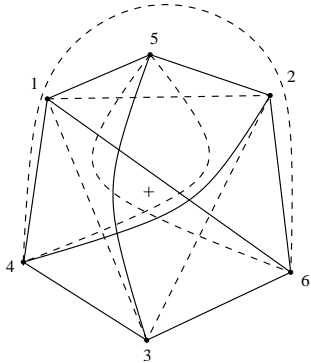} 
    \caption{$K_6$ with exactly one winding Birch partition.}  
    \label{fig-k6-one-wbp} 
\end{figure}

Motivated by recent work of Sch\"oneborn and Ziegler~\cite{schoeneborn05:_topol_tverb_theor},
and Remark~\ref{rem--bp-lower-bound} we have also studied the concept of
winding Birch partitions to obtain lower bounds
in the topological setting, see Hell~\cite{hell06:_tverb_fract_helly} for more details. The properties of Theorem~\ref{thm-number-birch-part} 
do not carry over to the topological setting. Hence a lower bound for the number of Tverberg partitions cannot be derived.
A computer project led to many examples of piecewise linear maps that have exactly one winding Birch partition for $k=2$; a smoothed version of one of them is
shown in Figure \ref{fig-k6-one-wbp}. There the only winding Birch
partition -- shown in broken lines -- is $\{1,2,3\}$ and $\{4,5,6\}$ with winding numbers $\pm
1$ resp.~$\pm 2$. For arbitrary dimension $d\geq 2$, 
note that any example for dimension $d$ can be
extended to an example in dimension $d+1$ using a construction
from Sch\"oneborn and Ziegler~\cite{schoeneborn05:_topol_tverb_theor}.

Let us end with two problems. Both are promising starting
points for future research.
\begin{problem}\label{prob-birch-gale}Relate the properties on
  the number $B_p(X)$ of Birch partitions to polytope theory.  Birch
  partitions show up while studying Gale diagrams; see Ziegler's
  textbook~\cite{ziegler95:_polytopes} for an introduction to Gale
  diagrams. In fact, a set $X$ of $k(d+1)$ points in $\R^d$ with
  $B_0(X)>0$ corresponds to a Gale diagram of a $k$-neighborly
  $(k-1)(d+1)$-dimensional simplicial polytope on $k(d+1)$ vertices.
\end{problem}
\begin{problem}\label{prob-birch-hc}
It is well-known that Radon's, Helly's, and
Carath\'eodory's theorem are closely related. Do the results on
the number of Birch partitions imply new Helly-type, or 
Carath\'eodory-type results?
\end{problem}

\section*{Acknowledgements}
Part of this work was stimulated by recent results on the colorful
simplicial depth by B\'ar\'any and Matou\v{s}ek~\cite{barany07:_quadr}, and Deza et.~al.~\cite{deza06:_colour}.
We thank Juliette Hell, Axel Werner, Carsten Schultz, G\"unter M.Ziegler,
and Rade \v{Z}ivaljevi\'c for helpful discussions and valuable remarks. All results are
part of my PhD thesis~\cite{hell06:_tverb_fract_helly}. 
This research was supported by the Deutsche Forschungsgemeinschaft 
within the European graduate program `Combinatorics, Geometry, and Computation' (No. GRK 588/2).



\end{document}